\documentclass[10pt,twocolumn,a4paper]{article} 
\usepackage{simpleConference}
\usepackage{times}
\usepackage{graphicx}
\usepackage{amssymb}
\usepackage{amsmath}
\usepackage{relsize}
\usepackage{url}

\usepackage{breakurl}
\usepackage[backref=page]{hyperref}

\newcommand{\conv}[2]{#1\star #2}
\newcommand{\norm}[1]{\left\lVert#1\right\rVert_2}
\newcommand{\abs}[1]{\left\lvert#1\right\rvert}
\newcommand{\R}[1]{\mathbb{R}^{#1}}
\newcommand{\CN}{\mathbb{C}^N}
\newcommand{\C}[1]{\mathbb{C}^{#1}}
\newcommand{\FN}{\mathcal{F}_{N}}
\newcommand{\Fr}[1]{\mathcal{F}_{#1}}
\newcommand*\conj[1]{\mkern 1.5mu\overline{\mkern-1.5mu#1\mkern-1.5mu}\mkern
1.5mu}
\newcommand{\layer}[1]{\mathrm{#1}}
\newcommand{\brac}[1]{\langle{#1}\rangle}

\newcommand{\bfb}{\mathbf{b}}
\newcommand{\bfc}{\mathbf{c}}
\newcommand{\bfx}{\mathbf{x}}
\newcommand{\bfX}{\mathbf{X}}
\newcommand{\bfy}{\mathbf{y}}
\newcommand{\bfw}{\mathbf{w}}
\newcommand{\bfz}{\mathbf{z}}
\newcommand{\bfZ}{\mathbf{Z}}
\newcommand{\bfu}{\mathbf{u}}

\begin{document}

\title{On the Equivalence of Convolutional and Hadamard Networks using DFT}

\author{Marcel Crasmaru \\
\today
\\
crasmarum@gmail.com  \\
}

\maketitle
\thispagestyle{empty}

\begin{abstract}
 In this paper we introduce activation functions that move the
 entire computation of Convolutional Networks into the frequency domain, where
 they are actually Hadamard Networks. To achieve this result we employ the
 properties of Discrete Fourier Transform. 
 
 We present some implementation details and experimental results, as well
 as some insights into why convolutional networks perform well in learning
 use cases.
\end{abstract}

\section{Introduction}

\begin{figure}[!h]
  \begin{center}
    \includegraphics[width=3.25in]{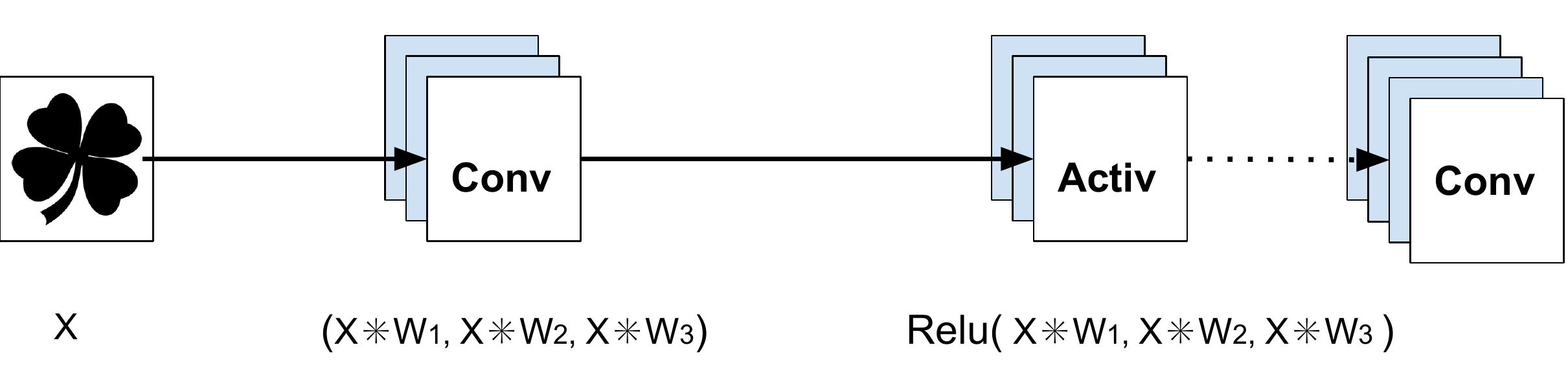}
  \end{center}

  \caption{\small Convolutional Network in the Space Domain.}
  \label{fig:conv}
\end{figure}

Convolutional networks use real numbers and operate in the spatial domain,
although their defining mathematical operation, that is, the convolution, is
easier to compute and understand after applying the Discrete Fourier Transform (\b{DFT}),
$\Fr{N}$, and moving the computation to the frequency domain. For example, as
you can see  in Figure ~\ref{fig:conv}, the first layer is convolving the 
input $X$ with three filters $W_1, W_2, W_3$. The output of the first layer 
$(X\star W_1, X\star W_2, X \star W_3)$ is then  forwarded to an activation layer, in this case
$\layer{Relu}.$

\begin{figure}[!h]
  \begin{center}
    \includegraphics[width=3.25in]{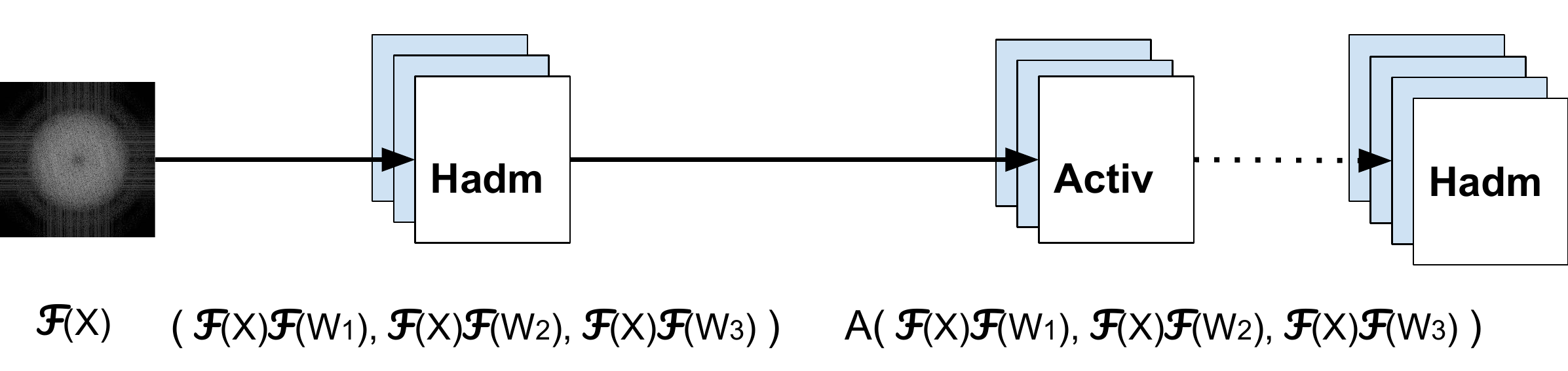}
  \end{center}

  \caption{\small Hadamard Network in the Frequency Domain.}
  \label{fig:hadm}
\end{figure}

In order to move the entire computation to the frequency domain we could apply
the DFT $\FN$ to the input $X$ and to the filters. We can then use the fact
that $\FN(X\star W) = \FN(X) \cdot \FN(W)$, that is, in the frequency domain
convolutions become the entrywise product of some complex numbers, which is known 
as the Hadamard product (see Figure~\ref{fig:hadm}).

In addition to simplifying the convolution calculations, there are many good
reasons for researching the convolutional network computation in the frequency
domain and use complex numbers for weights \cite{medium-article}. One compelling
reason is that real-valued networks cannot properly get close to the optimum on 
certain sets that are connected in
the complex space but not in the real one \cite{kronecker}.
Moreover, many researchers have observed that nets using complex values
for weights have a more robust and stable behavior in training 
\cite{unitary-rnn}. We could also mention some recent Computational Neuroscience
studies \cite{cortex}, that suggest that the entorhinal cortex might employ some 
sort of Fourier transform in its functionality.

The first goal of our paper is to investigate whether we can create
convolutional networks that operate entirely in the frequency domain and that
have all their layers ``commute" with the discrete Fourier transform in the
same way as convolutions. For example, one of the issues we have to address
is the fact that the currently used activation functions, for example,
$\layer{Relu}$, do not get mapped nicely to the frequency domain. What we would
like to find is a pair of easy-to-compute activation functions $A$ and $A'$
such that $\FN(A(\bfx))=A'(\FN(\bfx)).$ 

In Section \ref{sec:activation}, we will introduce a pair $\layer{A_N},
\layer{A_{\brac{N}}}$ of activation functions that commute with the DFT. We will
also show that all other widely used layers (fully connected, residual, 
batch normalization) are commutable with the DFT.

The second goal is to test that these convolutional networks converge for
some known standard data sets, for example, MNIST, and then see if we can use
Fourier analysis to get some insights on some hard open problems: why are
convolutional networks working well, is there a good reason why deep networks
generalize better than the shallow ones, etc. 

Experimental results show that a specific shallow network using
$\layer{A_{\brac{N}}}$ as an activation function converges fast in the frequency
domain:
we got $90\%$ accuracy on MNIST after $70$ epochs of Wirtinger gradient descent
back-propagation. This experimental result confirms the findings in \cite{adali}
and \cite{universal-apprx}, that is, as $\layer{A_{\brac{N}}}$ is bounded almost
everywhere but not holomorphic we should expect that feed-forward networks using this
activation function are capable of universally approximating any non-linear
complex mapping.

In this paper we show that any convolutional network built with layers that 
commute with the DFT is equivalent to a Hadamard
network. That is, the known technique of classifying “patterns” 
from input $\bfx$ using convolutional networks is equivalent to the technique of 
extracting the “recurring” information from the input $\bfx$ by applying DFT
$\FN(\bfx)$ then feeding it to a Hadamard network.

We propose that this equivalence might explain why convolutional networks work
well in areas such as image recognition. This equivalence also suggests that
predictive models using complex numbers for weights should perform beter than
the ones using only real numbers.
\section{Preliminaries}
\subsection{The Discrete Fourier Transform}

Let $N$ be an integer and suppose that $\bfx=(x_0,\ldots,x_{N-1})$ is an $N$
dimensional complex vector.
Let $\omega=\exp(\frac{-2\pi i}{N})$. Then the Discrete Fourier
Transform (DFT) $\mathcal{F}_N \colon \C{N} \to \C{N}$ is given by

\begin{equation} \label{def_dft}
X_k = \mathcal{F}_N(\bfx)_k=\sum_{j=0}^{N-1}x_j\omega^{jk}.
\end{equation}

The inverse DFT
$\mathcal{F}^{-1}_N \colon \C{N} \to \C{N}$ is then

\begin{equation*}
x_k = \mathcal{F}^{-1}_N(\bfX)_k={1\over N}\sum_{j=0}^{N-1} X_j \omega^{-jk}.
\end{equation*}

The DFT is a linear function, that is 
\begin{equation} \label{dft_liner}
\mathcal{F}_N(\alpha \bfx + \beta \bfy) = \alpha \mathcal{F}_N(\bfx) +
\beta \mathcal{F}_N(\bfy).
\end{equation}

Parseval Theorem states that if $\bfx, \bfw \in \C{N}$ then the following
equality holds:
\begin{equation}\label{parseval}
\sum_{j=0}^{N-1}x_j \conj{w_j} = \frac{1}{N} \sum_{j=0}^{N-1}
\mathcal{F}_N(\bfx)_j \conj{\mathcal{F}_N(\bfw)_j}.
\end{equation}

From this theorem one can easily prove the Plancherel Theorem that states that
if $\bfx,\bfX \in \CN$ and $\bfX=\FN(\bfx)$ then
\begin{equation}\label{plancheret}
\sum_{j=0}^{N-1}\abs{x_j}^2 = {1\over N}\sum_{j=0}^{N-1}\abs{X_j}^2
\end{equation}
and this is equivalent to 
\begin{equation}\label{Plancherel}
\norm{\bfx} = \frac{1}{\sqrt{N}} \norm{\FN(\bfx)}
\end{equation}
where $\norm{\bfx}=\sqrt{\sum_{j=0}^{N-1}x_j\conj{x_j}}.$

The circular convolution $\bfx \star \bfy \in \CN$ of two vectors $\bfx,\bfy \in
\CN$ is defined by
\begin{equation} \label{eqn:convo}
(\bfx \star \bfy)_k = \sum_{j=0}^{N-1} x_j y_{{k - j}\pmod N}.
\end{equation}

 The convolution property of the Fourier transform is
 \begin{equation} \label{eqn:conv_prop}
 \FN(\bfx \star \bfy) = \FN(\bfx)\cdot\FN(\bfy)
 \end{equation}
 where the multiplication indicated by the dot is element-wise.

\subsection{Wirtinger Derivatives and Gradients}

Given a function $g\colon \C{} \to \C{}$ differentiable in the real
sense, its Wirtinger derivative (equivalent to the areolar derivative
first introduced by Pompeiu in \cite{pompeiu})  at $c \in \C{}$ is
defined by:
\begin{equation*}
\frac{\partial g}{\partial z}(c) = 
\frac{1}{2}
\Big(\frac{\partial g}{\partial x}(c) - i \frac{\partial g}{\partial y}(c)\Big).
\end{equation*}

Its conjugate Wirtinger derivative is defined by:
\begin{equation*}
\frac{\partial g}{\partial \conj{z}}(c) = 
\frac{1}{2}
\Big(\frac{\partial g}{\partial x}(c) + i \frac{\partial g}{\partial y}(c)\Big).
\end{equation*}  

Note that a function $g(z) = u(x, y) + iv(x,y)$ is holomorphic (thus
satisfying the Cauchy Riemann equation $u_x=v_y, u_y=-v_x$)
if and only if its conjugate Wirtinger derivative vanishes, that is,
$\frac{\partial g}{\partial \conj{z}}(c) = 0$.

If we denote the conjugate of the complex number $z$  by $\conj{z}$ then the
following relations hold:
\begin{equation*}
\begin{split}
\conj{\frac{\partial g}{\partial z}(c)} = 
\frac{\partial \conj{g}}{\partial \conj{z}}(c)
\\\\
\conj{\frac{g}{\partial \conj{z}}(c)} = 
\frac{\partial \conj{g}}{\partial z}(c).
\end{split}
\end{equation*}  

For a function $f\colon \CN \to \C{}$ that is differentiable 
with respect to the real and imaginary parts of all of its inputs, the conjugate
gradient $\conj{\nabla}f(\bfc)$ is defined similarly to its real counterpart:
\begin{equation*} 
\conj{\nabla}f(\bfc)= \Big(\frac{\partial f}{\partial \conj{z_1}}(\bfc), \dots,
\frac{\partial f}{\partial \conj{z_N}}(\bfc)\Big).
\end{equation*} 

For any real-valued function (such as a cost function)
knowing its $\conj{\nabla}f(\bfc)$ gradient is a sufficient condition 
to minimize it, see for example \cite{wcalculus}.

The chain rules for two functions $f, g \colon \C{} \to \C{}$ that
are differentiable in the real sense, are:
\begin{equation*} 
\begin{aligned}
\frac{\partial f\circ g}{\partial z} = 
\Big(\frac{\partial f}{\partial z}\circ g\Big)\cdot\frac{\partial g}{\partial z}
+ \Big(\frac{\partial f}{\partial \conj{z}}\circ g\Big)
\cdot\frac{\partial \conj{g}}{\partial z}
\\
\frac{\partial f\circ g}{\partial \conj{z}} = 
\Big(\frac{\partial f}{\partial z}\circ g\Big)\cdot\frac{\partial
g}{\partial \conj{z}} + \Big(\frac{\partial f}{\partial \conj{z}}\circ g\Big)
\cdot\frac{\partial \conj{g}}{\partial \conj{z}}.
\end{aligned}
\end{equation*}

The above rules extend easily to functions of several variables. Given
two functions $f\colon \C{M} \to \C{}$ and $g\colon \C{N}
\to \C{M}$ differentiable in the real sense, the chain rules are:

\begin{equation} \label{def_chainrule}
\begin{aligned}
\frac{\partial f\circ g}{\partial z_k} = 
\sum_{m=1}^{M} \Big(\frac{\partial f}{\partial z_m}\circ
g \Big)\frac{\partial g_m}{\partial z_k} + 
\sum_{m=1}^{M}\Big(\frac{\partial f}{\partial \conj{z_m}}\circ 
g\Big)\frac{\partial \conj{g_m}}{\partial z_k}
\\
\frac{\partial f\circ g}{\partial \conj{z_k}} = 
\sum_{m=1}^{M} \Big(\frac{\partial f}{\partial z_m}\circ
g \Big)\frac{\partial g_m}{\partial \conj{z_k}} + 
\sum_{m=1}^{M}\Big(\frac{\partial f}{\partial \conj{z_m}}\circ 
g\Big)\frac{\partial \conj{g_m}}{\partial \conj{z_k}}
\end{aligned}
\end{equation}
where $k\in[1,N].$

\section{Convolutional Networks as Hadamard Networks in the Frequency
Domain}

For reasons we will provide in the next subsection and also for the sake
of simplicity, we will consider only complex valued networks operating on 
2D tensors. For example, a layer $\layer{F}$ with a $3$-channels tensor input
and a $2$-channels tensor output of channel size $256$ is a complex valued function 
$\layer{F}\colon \C{256\times 3} \to \C{256\times 2}.$ In general we will define
layers $\layer{G_{\brac{N}}}$ with $p$-channel tensor inputs and $q$-channel tensor
output of channel size $N$ as complex valued functions mapping 2D tensors, for example,
$\layer{G_{\brac{N}}} \colon \C{N \times p} \to \C{N \times q}.$

Notation-wise, we usually use upper indices for
channels and lower indices for vector/tensor components. We write a
$3$-channel tensor $\bfx\in \C{N\times 3}$ as $\bfx=(\bfx^1,\bfx^2,\bfx^3)$ and
by $x^2_j$ we mean the $j$ scalar component of vector $\bfx^2$, etc. We also
use the notation $\layer{G_N}$ for layers $\layer{G}$ used in the spatial domain
acting on tensors of channel size $N$ and $\layer{F_{\brac{N}}}$ for layers
$\layer{F}$ used in the frequency domain, etc.

Note also that Fourier transforms are always performed on the 1D components of
the 2D tensors: if $\bfx\in \C{N \times m}$ and $\bfx=(\bfx^1,\dots,\bfx^m)$
then
\begin{equation*}
\begin{split}
\FN(\bfx) &= (\FN(\bfx)^1,\ldots, \FN(\bfx)^m) \\ 
       &= (\FN(\bfx^1),\ldots,\FN(\bfx^m))
\end{split}
\end{equation*}
where, for each $1\leq j \leq m$ we compute $\FN(\bfx^j)$ as per definition
\eqref{def_dft}, etc.

\subsection{The Circular vs the Linear Convolution}
 Currently all convolutional layers used in the spatial domain use linear
 convolutions rather than circular ones. However, any linear convolution can
 be mapped to and computed from an 1D circular one. For example, given the 
 following 2D tensor $X$ and the $3\times3$ kernel matrix $W$:

\begin{equation*}
\begin{tabular}{ c c }
X =
  $\begin{bmatrix}
    \boldsymbol{x_0} & \boldsymbol{x_1} & x_2  \\
    \boldsymbol{x_3} & \boldsymbol{x_4} & x_5  \\
    x_6 & x_7 & x_8  \\
  \end{bmatrix}$

 & 

 W =
  $\begin{bmatrix}
    w_0 & w_1 & w_2  \\
    w_3 & \boldsymbol{w_4} & \boldsymbol{w_5}  \\
    w_6 & \boldsymbol{w_7} & \boldsymbol{w_8}  \\
  \end{bmatrix}$
 \end{tabular}
\end{equation*}
their 2D linear convolution of stride $1$ and same padding is the $3\times 3$
matrix $X \star W$ computed as $(X \star W)_{00} = x_0w_4 + x_1w_5 + x_3w_7 + x_4w_8$, \ldots, $(X \star W)_{11} =
\Sigma_{i=0}^{8}x_iw_i$, \ldots, etc. One can get the same result by padding
these matrices with zeros as shown here:
\[
\begin{tabular}{ c }

X =
  $\begin{bmatrix}
    0 & 0 & 0 & 0 & 0 \\
    0 & \boldsymbol{x_0} & \boldsymbol{x_1} & x_2 & 0 \\
    0 & \boldsymbol{x_3} & \boldsymbol{x_4} & x_5 & 0 \\
    0 & x_6 & x_7 & x_8 & 0 \\
    0 & 0 & 0 & 0 & 0 \\
  \end{bmatrix}$
\\ \\
  W =
  $\begin{bmatrix}
    w_0 & w_1 & w_2 & 0 & 0 \\
    w_3 & \boldsymbol{w_4} & \boldsymbol{w_5} & 0 & 0 \\
    w_6 & \boldsymbol{w_7} & \boldsymbol{w_8} & 0 & 0 \\
    0 & 0 & 0 & 0 & 0 \\
    0 & 0 & 0 & 0 & 0 \\
  \end{bmatrix}$

\end{tabular}
\]
and then flattening them in row major order to 1D tensors $\bfx, \bfw$:
\begin{alignat*}{16}
    & (0,\ && 0,\ && 0,\ && 0,\ && 0, \ && 0, \ && \boldsymbol{x_0}, \
    && \boldsymbol{x_1}, \ && x_2, \ && 0, \ && 0, \ && \boldsymbol{x_3}, \ &&
    \boldsymbol{x_4}, \ \dots ) \\
    & (w_0,\ && w_1,\ && w_2,\ && 0,\ && 0, \ && w_3, \ &&
    \boldsymbol{w_4}, \ && \boldsymbol{w_5}, \ && 0, \ && 0, \ && w_6, \ && \boldsymbol{w_7}, \ &&
    \boldsymbol{w_8}, \ \dots)
\end{alignat*}

One can observe that the dot product of the above tensors is precisely
$(X \star W)_{00}$. Moreover, one can see from ~\eqref{eqn:convo}, that this
value can be found by computing the circular convolution $\bfx \star \bfw^*$, where
$\bfw^*$ is obtained from the tensor $\bfw$ by performing $7$ left rotations and then
flipping the indices greater than $1$ left to right.

Note also that if $\bfx$ and $\bfy$ are two 1D tensors, if $\bfx$ has length
$|\bfx|$, and $\bfy$ has length $|\bfy|$, and both $\bfx$ and $\bfy$ are zero
padded to length $|\bfx| + |\bfy| - 1$, then their circular convolution matches
their linear convolution.

For these reasons henceforth we will only use 2D tensors and circular
convolutions, that is, the $\bfx \star \bfy$ notation will always mean the 
circular convolution. 

\subsection{The Convolution Layer}\label{seq_conv}
Given two $m$-channel tensors, that is, $\bfx, \bfw\in
\C{N\times m}$, we define their convolution as the sum of convolutions of their
1D components:
\begin{equation*}
\conv{\bfx}{\bfw} =\conv{(\bfx^1,\ldots,\bfx^m)}{(\bfw^1,\ldots,\bfw^m)} 
= \sum_{j=1}^m \conv{\bfx^j}{\bfw^j}.
\end{equation*}
By the linearity of DFT and the convolution property \eqref{eqn:conv_prop}, we
know $\FN(\conv{\bfx}{\bfw}) = \sum_{j=1}^m \FN(\bfx)^j\cdot\FN(\bfw)^j.$ 
Using this observation it follows that the convolutional layer
in the spatial domain:
\begin{equation*}
\layer{Conv_N}(\bfx,\bfw_1,\ldots,\bfw_n, \bfb) = (\conv{\bfx}{\bfw_1} +
\bfb^1,\ldots,\conv{\bfx}{\bfw_n} + \bfb^n)
\end{equation*}
commutes with the DFT:
\begin{equation*}
\begin{split}
\FN(&\layer{Conv_N}(\bfx,\bfw_1,\ldots,\bfw_n, \bfb)) \\
&=(\FN(\conv{\bfx}{\bfw_1}) + \FN(\bfb^1),\\
&\quad\quad\quad \ldots,\FN(\conv{\bfx}{\bfw_n}) + \FN(\bfb^n)) \\
&=\big(\sum_{j=1}^m \FN(\bfx)^j\cdot\FN(\bfw_1)^j + \FN(\bfb^1), \\
&\quad\quad\quad \ldots, \sum_{j=1}^m \FN(\bfx)^j\cdot\FN(\bfw_n)^j +
\FN(\bfb^n)\big),
\end{split}
\end{equation*}
where the input and the weights are $m$-channel tensors, that is,
$\bfx,\bfw_1,\ldots,\bfw_n \in \C{N\times m}$ and the bias $\bfb$ verifies that
$\bfb \in \C{N\times n}.$

Note that usually the tensors used as weights in the convolutional layer
have low dimensionality and they are called ``kernels'', for example $3\times 3$
matrices are widely used etc. As we have observed in the previous subsection,
computing linear convolutions with these kernels is equivalent to performing 1D
circular convolutions, after some suitable zero padding. 
Based on this observation we define a $Pad$ function that takes a linear
pattern, that is, an ordered set of indices $K\subseteq \{0,\ldots,N-1\}$, and
$k = |K|$ tensor weights to create (via some zero padding) a length $N\times m$
tensor of weights.
\newline
More precisely we define $Pad_N \colon \C{k\times m}\times \mathcal{P}(\{0,
N-1\}) \to \C{N\times m} $ by
\begin{equation*}
  Pad_N(\bfw, K)_j^l = 
  \begin{cases} 
      w_i^l & \text{ if } {K[i] = j} \\\\
      0 & \text{ otherwise}
   \end{cases}
\end{equation*}
for $0 \leq j < N$ and $1 \leq l \leq m$.

For example, given a $3\times 3$ kernel of weights $\bfw \in \C{9\times 1}$ and 
$K=(0, 1, 2, 5, 6, 7, 10, 11, 12)$ then $Pad_{25}(\bfw, K)$ is equal to
\begin{alignat*}{16}
    & (w_0,\ && w_1,\ && w_2,\ && 0,\ && 0, \ && w_3, \ &&
    w_4, \ && w_5, \ && 0, \ && 0, \ && w_6, \ && w_7, \ &&
    w_8, \ && 0,\dots).
\end{alignat*}

We can then extend convolutional layers in the spatial domain having a kernel
$K\subseteq \{0,\ldots,N-1\}$ as follows.
\begin{equation*}
\begin{split}
&\layer{Conv_N} (\bfx,K, \bfw_1,\ldots,\bfw_n, \bfb) = \\
 &\quad\quad\layer{Conv_N} (\bfx,Pad_N(\bfw_1, K),\ldots,Pad_N(\bfw_n, K),
 \bfb).
\end{split}
\end{equation*}

\subsection{The Hadamard Layer}

The Hadamard layer $\layer{H_{\brac{N}}}$ corresponding to
the convolutional layer $\layer{Conv_N}$ is simply defined
so that the DFT carries over from the spatial domain to the frequency
domain:
$\Fr{N}(\layer{Conv_N}(\bfz,K,
\bfw_1,\ldots,\bfw_q,\bfb))=\layer{H_{\brac{N}}}(\Fr{N}(\bfz),K,
\bfw_1,\ldots,\bfw_q,\Fr{N}(\bfb)).$
\newline
Formally, a Hadamard layer with a $p$-channels input tensor
$\bfz=(\bfz^1,\dots,\bfz^{p})\in \C{N\times p}$ and outputting a $q$-channels
tensor $\bfy=(\bfy^1,\ldots,\bfy^q)\in \C{N\times q}$ by using the kernel
$K\subseteq \{0, N-1\}$, the weights 
$\bfw_1,\dots, \bfw_{q} \in \C{|K|\times p}$ and the biases $\bfb \in
\C{N\times q}$ is defined as the complex multivariate function
$\layer{H_{\brac{N}}} \colon \C{N\times p} \times \mathcal{P}(\{0,
N-1\}) \times \C{|K|\times p} \times
\C{N\times q} \to \C{N\times q} $ given by
\begin{equation}
\bfy^j = \sum_{i=1}^p \bfz^i \cdot \FN(Pad(\bfw_j^i,K)) + \bfb^j
\end{equation}
where $1\leq j \leq q$, ``$\mathbf{\cdot}$'' is the Hadamard product of two
vectors, $\FN$ is the discrete Fourier transform and $Pad$ is the function
defined previously.

The Wirtinger derivatives of this layer with respect to input and
weights can be computed by applying the chain rule \eqref{def_chainrule}, and
the definition of the discrete Fourier transform \eqref{def_dft}. More precisely,
assume we have back-propagated the Wirtinger gradients for the layer $\layer{G}$ 
where $L \colon \C{N\times q} \to \R{}$ and $L = \layer{G} \circ
\layer{H_{\brac{N}}}$ is the loss function. The back-propagation rules for the
layer $\layer{H_{\brac{N}}}$ and the input variable $z^i_k$ are:
\begin{equation*}
\frac{\partial \layer{G}\circ\layer{H_{\brac{N}}}}{\partial z^i_k} =
\sum_{j=1}^q \FN(Proj(\bfw_j^i,K))_k \frac{\partial \layer{G}}{\partial
z^j_k}\circ
\layer{H_{\brac{N}}}
\end{equation*}
and
\begin{equation*}
\frac{\partial \layer{G}\circ\layer{H_{\brac{N}}}}{\partial \conj{z^i_k}} =
\sum_{j=1}^q
\conj{\FN(Proj(\bfw_j^i,K))_k} \frac{\partial \layer{G}}{\partial
\conj{z^j_k}}\circ
\layer{H_{\brac{N}}}
\end{equation*}
where $1\leq i \leq p$ and $0\leq k < N$.

\subsection{The Input Layer}
In order to move the entire computation in the frequency domain we will simply
apply the Fourier transform to the input tensors. That is, the input tensor
$\bfz = (\bfz^1, \dots, \bfz^k) \in \C{N\times k}$ is mapped to
$\bfZ = \FN(\bfz)=(\FN(\bfz^1),\dots,\FN(\bfz^k)).$ For example, a CIFAR-10
input tensor $\bfx=(\bfx^1, \bfx^2, \bfx^3) \in \R{1024 \times 3}$ is mapped to
$\bfX=(\Fr{1024}(\bfx^1),\Fr{1024}(\bfx^2), \Fr{1024}(\bfx^3))$, etc. 
\linebreak 
Formally we define
the layer $\layer{In_{\brac{N}}}$ to be the function $\layer{In_{\brac{N}}} \colon \C{N\times k}
\to \C{N\times k}$ given by
\begin{equation}
\layer{In_{\brac{N}}}(\bfx) = (\FN(\bfx^1),\dots,\FN(\bfx^k)).
\end{equation}
Note that by using fast Fourier transform algorithms we can compute
$\layer{In_{\brac{N}}}(\bfx)$ in $O(kN\log{N})$ time. Given the training/testing
examples/data points $(\bfx_{1}, y_{1}),\dots,(\bfx_{m}, y_{m})$ a faster
approach is to pre-compute them ahead of time: 
$(\FN(\bfx_{1}), y_{1}),\dots,(\FN(\bfx_{m}), y_{m}).$

\subsection{The Activation Layer}\label{sec:activation}

Let $\layer{A_N}$ be the function $\layer{A_N} \colon \CN \to \CN$ given by
\begin{equation}
\layer{A_N}(\bfx) = \begin{cases} 
      \dfrac{\bfx}{\sqrt{N} \norm{\bfx}} & x\ne \boldsymbol{0} \\
      0 & x=\boldsymbol{0} .
   \end{cases}
\end{equation}

One can easily observe that $\layer{A_N}$ is bounded for all \mbox{$\bfx
\in\CN$}:
\begin{equation*}
\norm{\layer{A_N}(\bfx)} \le \frac{1}{\sqrt{N}}.
\end{equation*}

Moreover, for a given function $g \colon \CN \to \CN$ and from equation
\eqref{Plancherel}, we know:
\begin{equation*}
\begin{split}
\FN(\layer{A_N}(g(\bfx)))_k & =\sum_{j=0}^{N-1}\layer{A_N}
(g(\bfx))_j\omega^{jk}
\\
 & =\sum_{j=0}^{N-1}\frac{g(\bfx)_j}{\sqrt{N} \norm{g(\bfx)}}\omega^{jk} \\
 & = \frac{1}{\sqrt{N} \norm{g(\bfx)}} \sum_{j=0}^{N-1}g(\bfx)_j\omega^{jk} \\
 & = \frac{\FN(g(\bfx))_k}{\norm{\FN(g(\bfx))}}.
\end{split}
\end{equation*}

This means that for any function $g \colon \CN \to \CN$ the
following equation holds:
\begin{equation*}
\FN(\layer{A_N}(g(\bfx))) = \begin{cases} 
      \dfrac{\FN(g(\bfx))}{\norm{\FN(g(\bfx))}} & g(\bfx)\ne \boldsymbol{0} \\
      0 & g(\bfx)=\boldsymbol{0}.
   \end{cases}
\end{equation*}
\newline
Based on the above observation we define the activation layer in the frequency
domain $\layer{A_{\brac{N}}} \colon \C{N} \to \C{N}$ to~be
\begin{equation}
\layer{A_{\brac{N}}}(\bfz) = \begin{cases} 
      \dfrac{\bfz}{\norm{\bfz}} & \bfz\ne \boldsymbol{0} \\
      0 & \bfz=\boldsymbol{0} 
   \end{cases}
\end{equation}
and we know $\layer{A_N}$ and $\layer{A_{\brac{N}}}$ commute with the DFT:
\begin{equation*}
\FN(\layer{A_N}(\bfz))=\layer{A_{\brac{N}}}(\FN(\bfz)).
\end{equation*}

Let
$\layer{A_{\brac{N}}}(\bfz)=(f_0(\bfz),\dots,f_{N-1}(\bfz))$
and $\bfz=(z_0,\ldots,z_{N-1})$. For each $k\in[0,N-1]$ we know 
$f_k\colon \CN \to \C{}$ and  $z_k = x_k +iy_k$.
The partial derivatives of $f_k(\bfz)$ with respect to the variables
$x_k$ and $y_k$ are:
\begin{equation*} 
\begin{aligned}
&\dfrac{\partial f_k}{\partial x_k}(\bfz)=\dfrac{\partial}{\partial
x_k}\dfrac{x_k + i y_k}{\norm{\bfz}} = \dfrac{\norm{\bfz}^2 -
x_k^2}{\norm{\bfz}^3} - i\dfrac{x_ky_k}{\norm{\bfz}^3}
\\
&\dfrac{\partial f_k}{\partial y_k}(\bfz)=\dfrac{\partial}{\partial
y_k}\dfrac{x_k + i y_k}{\norm{\bfz}} = -\dfrac{x_ky_k}{\norm{\bfz}^3} 
+ i\dfrac{\norm{\bfz}^2 - y_k^2}{\norm{\bfz}^3}.
\end{aligned}
\end{equation*} 
The Wirtinger derivatives for $f_k(\bfz)$ follow from the above
equations and their definitions:
\begin{equation*} 
\begin{aligned}
&\dfrac{\partial f_k}{\partial z_k}(\bfz) =
\dfrac{1}{\norm{\bfz}} - \frac{1}{2}\dfrac{|z_k|^2}{\norm{\bfz}^3}
\\
&\dfrac{\partial f_k}{\partial \conj{z_k}}(\bfz) =
-\frac{1}{2}\dfrac{z_k^2}{\norm{\bfz}^3}
\\
&\dfrac{\partial f_j}{\partial z_k}(\bfz) =
-\frac{1}{2}\dfrac{z_j \conj{z_k}}{\norm{\bfz}^3} \text{, when }j\neq k
\\
&\dfrac{\partial f_j}{\partial \conj{z_k}}(\bfz) =
-\frac{1}{2}\dfrac{z_j z_k}{\norm{\bfz}^3} \text{, when }j\neq k.
\end{aligned}
\end{equation*} 
\newline
Note that the $\layer{A_{\brac{N}}}$ activation function is bounded and
differentiable in the real sense but it is not differentiable in the complex
sense as its $\dfrac{\partial f_j}{\partial \conj{z_k}}$
derivatives do not vanish, hence $\layer{A_{\brac{N}}}$ is not holomorphic. 
\newline
Note also that Georgiou et al introduced the following activation function in
\cite{afunction}:

\begin{equation*} 
f(z) = \frac{z}{c + \frac{1}{r}|z|}
\end{equation*}
where $c$ and $r$ are real positive constants. This function maps a point $z$
on the complex plane to a unique point $f(z)$ on the open disc \
$\{ z \colon |z| < r\}$. One way to extend this activation function to $\CN$,
while making it commute with the DFT, is the following:
\begin{equation*}
\layer{A_{N, c, r}}(\bfz) = 
      \dfrac{\bfz}{c + \frac{\sqrt{N}}{r} \norm{\bfz}}
\end{equation*}
and we can similarly show that $\layer{A}_{N, c, r}$ has the
property that
\begin{equation*}
\FN(\layer{A_{N, c, r}}(\bfz)) = 
      \dfrac{\FN(\bfz)}{c + \frac{1}{r} \norm{\FN(\bfz)}}.
\end{equation*}

Henceforth we will extend these activation functions to 2D tensors in the usual
way, that is, if $\bfz= (\bfz^1,\ldots,\bfz^m)\in\C{N\times m}$  then
\begin{equation}
\layer{A_{\brac{N}}}(\bfz) =
(\layer{A_{\brac{N}}}(\bfz^1),\ldots,\layer{A_{\brac{N}}}(\bfz^m)).
\end{equation}

\subsection{The Fully Connected Layer}

Parseval Theorem states that if $\bfz, \bfw \in \C{N}$ then the following
equality holds:
\begin{equation*}
\sum_{j=0}^{N-1}z_j \conj{w_j} = \frac{1}{N} \sum_{j=0}^{N-1}
\mathcal{F}_N(\bfz)_j \conj{\mathcal{F}_N(\bfw)_j}
\end{equation*}
that is, the complex dot product $\bfz \bullet \bfw$ equals $\frac{1}{N}
\mathcal{F}_N(\bfz) \bullet \mathcal{F}_N(\bfw).$ 
Based on this equality it makes sense to define Fully
Connected layers $\layer{FC_{\brac{N}}}$ in the frequency domain as follows:
\begin{equation}
\begin{aligned}
&\layer{FC_{\brac{N}}}(\bfz, \bfw_1,\dots, \bfw_k) =
\\ &\Big(\frac{1}{N} \sum_{j=0}^{mN-1} z_j \conj{w_{1j}},
\dots,
\frac{1}{N} \sum_{j=0}^{mN-1} z_j \conj{w_{kj}} \Big)
\end{aligned}
\end{equation}
where $\layer{FC_{\brac{N}}} \colon \C{N\times m} \times \C{N\times m\times k} \to
\C{k}$, $\bfz$ is an $m$-channels input and $\bfw_1,\dots, \bfw_k \in \C{N\times
m}$ are the fully connected layer's filters/weights.

This definition has the property that the well-known fully connected layers
operating in the spatial domain are mapped naturally to the frequency domain.
That is, if $\bfx = (\bfx^1,\ldots,\bfx^m)\in \R{N\times m}$ is an $m$-channels
input tensor and $\bfw_1,\dots, \bfw_k \in \R{N\times m}$ are the weights of a
standard fully connected layer $\layer{FC_N}$, then the layer's output verifies
the following relation:
\begin{equation*}
\begin{aligned}
\layer{FC_N}(\bfx,&\bfw_1,\ldots,\bfw_k)
\\
=&\Big(\sum_{j=0}^{mN-1} x_j w_{1_j}, \dots, \sum_{j=0}^{mN-1} x_j w_{k_j} \Big)
\\
=&\Big(\sum_{j=0}^{mN-1} x_j \conj{w_{1_j}},\dots, \sum_{j=0}^{mN-1} x_j
\conj{w_{k_j}} \Big)
\\
=&\Big(\sum_{i=1}^{m} \bfx^i \bullet \bfw^i_1,\dots, \sum_{i=1}^{m} \bfx^i
\bullet \bfw^i_k \Big)
\\
=&\frac{1}{N}\Big(\sum_{i=1}^{m} \FN(\bfx^i) \bullet \FN(\bfw^i_1),
\\ &\quad\quad\quad\quad \dots, \sum_{i=1}^{m} \FN(\bfx^i) \bullet
\FN(\bfw^i_k)\Big)
\\
=& \layer{FC_{\brac{N}}}(\FN(\bfx),\FN(\bfw_1),\dots, \FN(\bfw_k))
\\
=& \layer{FC_{\brac{N}}}(\FN(\bfx, \bfw_1,\dots, \bfw_k)).
\end{aligned}
\end{equation*}

Note that formally we define $\layer{FC_N}$ as a complex valued function 
$\layer{FC_N} \colon \C{N\times m} \times \C{N\times m\times k} \to
\C{k}$ and 
\begin{equation*}
\layer{FC_N}(\bfx,\bfw_1,\ldots,\bfw_k)=\Big(\sum_{i=1}^{m} \bfx^i \bullet
\bfw^i_1,\dots, \sum_{i=1}^{m} \bfx^i \bullet \bfw^i_k \Big)
\end{equation*}
and we have that
\begin{equation*}
\layer{FC_N}(\bfx,\bfw_1,\ldots,\bfw_k)
=\layer{FC_{\brac{N}}}(\FN(\bfx,\bfw_1,\ldots,\bfw_k)).
\end{equation*}

The Wirtinger gradients of the $\layer{FC_{\brac{N}}}$ layer 
 can be computed easily at back-propagation time by applying the chain rule
\eqref{def_chainrule} and the trivial equalities $\partial
(z\conj{w}) / \partial \conj{w} = z$ and $\partial
(z\conj{w})/{\partial z} = \conj{w}$.

\subsection{The Output Layer and the Cross Entropy Loss} \label{sec:loss}

Let $\layer{Out_N}$ be the function $\layer{Out_N} \colon \CN \to \CN$ given
by
\begin{equation}
\layer{Out_N}(\bfz) = \begin{cases} 
      \dfrac{\bfz}{\norm{\bfz}} & \bfz\ne \boldsymbol{0} \\
      \\
      (\sqrt{1/N},\dots,\sqrt{1/N}) & \bfz=\boldsymbol{0} .
   \end{cases}
\end{equation}

For all $\bfz = (z_{0}, \dots,z_{N-1})\in \CN$ we have that $1 =
\mathlarger{‎‎\sum}_{i = 0}^{N-1} |\layer{Out_N}(\bfz)_j|^2.$ Thus the
probability of a particular outcome $j \in 0,\dots,N-1$ is 
\begin{equation*}
P(y = j \mid \bfz) =
|\layer{Out_N}(\bfz)_j|^2 = z_j \cdot \conj{z_j} / \norm{\bfz}^2.
\end{equation*}
That is, the squared absolute value of the $\layer{Out_N}(\bfz)_j$ amplitude. 
Incidentally, this way of defining probabilities is similar to the Born Rule
in Quantum Mechanics: the probability of a particular outcome is the squared
absolute value of a certain amplitude.
\newline 
Given a pair $(\bfz^0, \bfy^0) \in
\CN \times \R{N}_{+}$ with $\sum_{j=0}^{N-1}y^0_j=1$, the cross entropy loss $L
\colon \CN \times \R{N}_{+} \to \CN$ is
\[
L(\bfz^0, \bfy^0)=\sum_{k=0}^{N-1} -y^0_k\log(z^0_k \conj{z^0_k} /
\norm{\bfz^0}^2).
\]

We can show that the cost functions $C_k \colon \CN \to \R{}$, where $0\leq k <
N$, given by $C_k(\bfz)=\log(z_k \conj{z_k} / \norm{\bfz}^2)$ have the
following Wirtinger derivatives
\begin{equation*}
\begin{aligned}
&\dfrac{\partial C_k}{\partial z_k}(\bfz) = \conj{z_k}(1/|z_k|^2
-1/\norm{\bfz}^2)
\\
&\dfrac{\partial C_k}{\partial \conj{z_k}}(\bfz) = z_k(1/|z_k|^2
-1/\norm{\bfz}^2)
\\
&\dfrac{\partial C_k}{\partial z_j}(\bfz) = -\conj{z_j} / \norm{\bfz}^2, \text{
when } k\neq j
\\
&\dfrac{\partial C_k}{\partial \conj{z_j}}(\bfz) = -z_j / \norm{\bfz}^2,  \text{
when } k\neq j.
\end{aligned}
\end{equation*}
The above equalities let us compute $\conj{\nabla}(L)$, the Wirtinger
conjugate gradient of the cross entropy loss function, $L$. As $L$ is a
real-valued function we can then use gradient descent algorithms to minimize it, 
see \cite{wcalculus}, for more information.

Note that 
\begin{equation*}
\layer{Out_N}=\layer{Out_{\brac{N}}},
\end{equation*}
that is, we will use the same output layer
in both the spatial and the frequency domains rather than having
$\layer{Out_{\brac{N}}}\equiv (1 /\sqrt{N}) \FN(\layer{Out_N})$. The reason for
this is that while the $\layer{Out_N}$ layer commutes with the DFT, the 
Hirschman entropic uncertainty principle tells us that we cannot minimize $L$, 
the cross entropy loss, in both domains at the same time. 

More precisely, let $\bfz\in \CN$ and for $0\leq j < N$ let $p_j =
|\layer{Out_N}(\bfz)_j|^2$ and $q_j =
\frac{1}{N} |\FN(\layer{Out_N}(\bfz))_j|^2$. We then know $\sum_{j =
0}^{N-1} p_j= 1$ and by Parseval Theorem we also know $\sum_{j =
0}^{N-1} q_j = 1$, however the Hirschman uncertainty principle tells us that
\begin{equation*}
 -\mathlarger{‎‎\sum}_{j = 0}^{N-1} p_j\log(p_j)  -  \mathlarger{‎‎\sum}_{j =
 0}^{N-1}q_j\log(q_j) \geq \log(N).
\end{equation*}
Roughly speaking, we cannot minimize~$L$, the cross entropy loss, in both the
space and the frequency domains, at the same time.

\subsection{Putting it All Together}

Let's start with an example by defining a minimal convolutional network $\layer{Net}$ in
the space domain in order to create a model for the MNIST dataset. Let this
$\layer{Net}$ consist of a convolution layer, say with $64$ filters of
kernel $K=\{0, 1, 2, 28, 29, 30, 56, 57, 58\}$,
an activation layer $\layer{A_{784}}$ and a fully connected layer $\layer{FC}$:
\begin{equation*} 
\begin{aligned}
\layer{Net}(&\bfx,K,\bfw_1\ldots,\bfw_{64},\bfb,\bfu) =
\\
&\layer{Out_{10}}(\layer{FC_{784}}(\layer{A_{784}}(\conv{\bfx}{Pad(K, \bfw_1)} +
\bfb^1,
\\
&\ldots, \conv{\bfx}{Pad(K, \bfw_{64})} + \bfb^{64})), \bfu^1,\ldots,\bfu^{10}))
\end{aligned}
\end{equation*}
where $\bfx\in \R{784}$, $\bfw_1,\ldots,\bfw_{64}\in \C{9}$,
$\bfb\in \C{784 \times 64}$ and $\bfu\in \C{784 \times 64 \times 10}$.

As we have already observed in the previous subsection, all $\layer{Net}$'s
layers are commuting with the DFT, therefore we can write
\begin{equation*} 
\begin{split}
\layer{Net}
&=\layer{Out_{\brac{10}}}(\layer{FC_{784}}(\layer{A_{784}}(\ldots)))
\\&=\layer{Out_{\brac{10}}}(\layer{FC_{\brac{784}}}(\Fr{784}(\layer{A_{784}}(\ldots))))
\\&\ldots
\\&=\layer{Out_{\brac{10}}}(\layer{FC_{\brac{784}}}(\layer{A_{\brac{784}}}
(\layer{H_{\brac{784}}}(\layer{In_{\brac{784}}}))))
\\&=
\layer{Out_{\brac{10}}}\circ\layer{FC_{\brac{784}}}\circ\layer{A_{\brac{784}}}
\circ\layer{H_{\brac{784}}}\circ\layer{In_{\brac{784}}}.
\end{split}
\end{equation*} 
That is, we can perform all the computation in the frequency domain, where
convolutions are mapped to Hadamard products which are simpler to
deal with. Moreover, we can minimize the cross entropy loss
by back-propagating in the frequency domain the Wirtinger gradients for inputs,
weights and biases. Another observation is that the only requirement for the
final $\layer{Out_N}=\layer{Out_{\brac{N}}}$ layer is to be a map from a complex
valued vector to a probability vector.

The above result is easily generalized to any convolutional network built with
layers that commute with the DFT. Therefore, we conclude that any convolutional 
network built with layers that commute with the Discrete Fourier Treansform is
equivalent to a Hadamard Network.

First consequence of this equivalence is that predictive models using complex
numbers for weights should perform better as the space of solutions is larger.

Second consequence is that performing the computation in the frequency
domain should be faster as convolutions become Hadamard products and we can
precompute the DFT of the inputs. 

Third consequence is that the above convolutional
networks are equivalent to multivariate rational expressions in the frequency
domain, and learning via gradient descent is a technique for
interpolation.

Last but lot the least, the technique we present of extracting the
``recurring'' information in the input $\bfx$ by applying DFT $\FN(\bfx)$ then
feeding it in a Hadamard network is equivalent to the known technique of 
extracting ``patterns'' from input $\bfx$ via convolutional networks.

\section{Implementation and Experimental Results}
In order to properly perform the back-propagation of gradients in the frequency
domain, we need support for complex number operations and also a way to keep
track of both the Wirtinger and the Wirtinger conjugate gradients for all inputs. 
For weights and biases, just keeping track of the Wirtinger
conjugate gradient is enough. The current machine learning frameworks
(Tensorflow, Torch, etc.) have limited support for both complex
differentiation and computing Fourier Transforms on GPUs, therefore the author
has run experiments on custom software. The hardware used was based on a NVIDIA
GTX-1070-Ti GPU.

\subsection{MNIST}

We have obtained $90\%$ accuracy on MNIST after $70$ epochs using a mini-batch
Wirtinger gradient descent of batch size $100$ and a shallow network in the
frequency domain. This network has only one Hadamard layer with $50$ out
channels and a kernel $K$ of size $7\times 7$ with $K=\{0,1,2,3,4,5,6,
28,29,30,31,32,33,34,\cdots,168(=28*6),169,170,171,172,173,174\}$.
This Hadamard layer is directly connected to an activation layer, followed by a
$10$ out channels fully connected layer, followed by an output layer:
\begin{equation*}
\begin{aligned}
\layer{In_{\brac{784}}} &\to \layer{H_{\brac{784, 50, 7 \times 7}}} \to
\layer{A_{\brac{784}}}
\\
&\to  \layer{FC_{\brac{784, 10}}} \to \layer{Out_{\brac{10}}}.
\end{aligned}
\end{equation*}
The above net has $394510$ complex parameters that are initialized using
the circularly symmetric complex Gaussian distribution, see Section
\ref{sec:weight-init} in the Annex for details.

We slightly improved the accuracy to $91.1\%$ by increasing the number of
Hadamard filters to $100$ and training the net for $200$ epochs.

\section{Conclusions and Further Work}

In this paper, we first investigated whether we can modify
convolutional networks so that the entire computation moves from the space
domain into the frequency domain and from convolutions to Hadamard products.

We found out that there are indeed certain activation functions and additional 
layers that commute with the DFT and make computation in the frequency domain
possible. 

Next, we did some experiments and found that there are shallow
Hadamard networks in the frequency domain that converge fast to a model for the
MNIST data.

Finally, we showed that any convolutional network built with layers that commute 
with the DFT is equivalent to a Hadamard Network. This result suggests that
using complex numbers for weights creates better predictive models. Another
consequence of this equivalence is that these convolutional networks are
multivariate rational expressions in the frequency domain, and learning via 
gradient descent is a technique for interpolation.

This paper opens some interesing lines of research and here are some questions
that might deserve some further theoretical and experimental investigation:

1. Could we build state of the art Hadamard predictive models for larger data
sets, say CIFAR10? Can we do it with a shallow Hadamard network? Do we need data
augmentation?

2. If $\bfb=(c_0,c_0,\ldots,c_0)\in\CN$ then $\FN(\bfb)=(N\cdot
c_0,0,\ldots,0),$ and $\FN(\bfx\star \bfw + \bfb)=\FN(\bfx)\cdot \FN(\bfw) +
(N\cdot c_0,0,\ldots,0).$ Does this mean that using a bias in the convolution 
layer is unnecessary as it only acts as a high pass filter for the stationary
frequency? Or should the bias be multivariate?

3. Is the divider layer proposed in the Section \ref{sec:divider} playing
the same role for Hadamard networks as the drop out layer is playing for
convolutional networks?

4. Batch normalization and residual layers are both
commuting with the DFT, see the Sections \ref{sec:batch-norm} and \ref{residual}
in the Appendix. Are these layers speeding up the learning in deep Hadamard 
networks?

5. Let $c\in\R{+}$ be a positive constant. The cost functions~$C_k$
and their Wirtinger gradients introduced in the Section \ref{sec:loss} are
Lipschitz  on the domain $\{\bfz : \norm{\bfz}>c\}.$ Could we use this
observation to analyse the gradient descent convergence?

\bibliographystyle{abbrv}

\section{Appendix}

\subsection{A Dropout Equivalent: the Divider Layer} \label{sec:divider}

If $N$ is even and $\bfx\in \C{N}$ then let $\bfx_{E} = (x_0, x_2, x_4,
\dots, x_{N-1})$ and $\bfx_{O} = (x_1, x_3, x_5, \dots, x_{N-2}).$ It is a
well known fact that the following recurrences hold for all $0 \leq k < N/2$:
 \begin{equation*}
 \begin{aligned}
 {\Fr{N}(\bfx)}_k &= {\Fr{N/2}(\bfx_{E})}_k + \omega^{k} {\Fr{N/2}(\bfx_{O})}_k 
 \\
 {\Fr{N}(\bfx)}_{k + N/2} &= {\Fr{N/2}(\bfx_{E})}_k - \omega^k
 {\Fr{N/2}(\bfx_{O})}_k .
 \end{aligned}
 \end{equation*}
 The above equalities are equivalent to the following relations:
  \begin{equation*}
 \begin{aligned}
 {\Fr{N/2}(\bfx_{E})}_k &= \frac{1}{2}  ({\Fr{N}(\bfx)}_k + {\Fr{N}(\bfx)}_{k +
 N/2})
 \\
 {\Fr{N/2}(\bfx_{O})}_k &= \frac{\omega^{-k}}{2} ({\Fr{N}(\bfx)}_k -
 {\Fr{N}(\bfx)}_{k + N/2}) .
 \end{aligned}
 \end{equation*}
 
Then, we can define the Divider layer as the function $\layer{Div_{\brac{2N}}}
\colon \C{2N} \to \C{N}\times\C{N}$
\begin{equation*}
  \layer{Div_{\brac{2N}}}(\bfx)_k = \begin{cases} 
      \frac{1}{2} (x_k + x_{k + N/2}) & k < N \\\\
      \frac{\omega^{-k}}{2}(x_k - x_{k + N/2}) & N \leq k < 2N 
   \end{cases}
\end{equation*}
having the following property for all $\bfx \in \C{2N}$:
\begin{equation}
  \layer{Div_{\brac{2N}}}(\Fr{2N}(\bfx)) = (\Fr{N}(\bfx_E), \Fr{N}(\bfx_O)).
\end{equation}
In words, the Divider layer halves the dimension of the channels by doubling
their numbers while preserving the data under the discrete Fourier transform.

\subsection{The Batch Normalization Layer Commutes with DFT}
\label{sec:batch-norm}

Rather than considering batch normalization in the complex domain to be
equivalent to whitening 2D vectors, as Trabelsi et al have proposed in
\cite{deep-cnn}, we would like to analyse how the 
batch normalization layer is mapped into the frequency domain under the
discrete Fourier transform. 

Let $\bfz_1,\ldots,\bfz_m$ be a batch of $m$ input tensors with $p$ channels /
feature layers such that $\bfz_j=(\bfz_j^1,\ldots,\bfz_j^{p})\in \C{N\times p}$.
For each layer $1 \leq k \leq p$ let the mean $\mu_k$ be
\begin{equation*}
\mu_k = \frac{1}{m} \sum_{j=1}^{m}\bfz_j^k.
\end{equation*}
Obviously, $\mu_k\in \C{N}$ for all $1\leq k \leq p.$ Similarly, let the
variance $\sigma_k^2$~be given by formula
\begin{equation*}
\sigma_k^2 = \frac{1}{mN} \sum_{j=1}^{m} (\bfz_j^k - \mu_k)\bullet(\bfz_j^k -
\mu_k)
\end{equation*}
where ``$\bullet$'' is the inner product of complex vectors etc.
Note that the variance $\sigma_k^2$ is a real number.
By using the DFT linearity \eqref{dft_liner} and the Parseval Theorem
\eqref{parseval} one can show that
\begin{equation}
\begin{aligned}
&\FN(\mu_k) = \frac{1}{m} \sum_{j=1}^{m} \FN(\bfz_j^k), \text{ and}
\\
&\sigma_k^2 = \frac{1}{mN^2} \sum_{j=1}^{m} (\FN(\bfz_j^k) - \FN(\mu_k))
\\ & \quad\quad\bullet(\FN(\bfz_j^k) - \FN(\mu_k)) \label{sigma_fft}.
\end{aligned}
\end{equation}
In other words, the mean $\mu_k^\prime$ of $\FN(\bfz_1^k),\ldots,\FN(\bfz_m^k)$
 is equal to $\FN(\mu_k)$ while their variance ${\sigma_k^2}^\prime$ is  $N$
times $\sigma_k^2.$ The $1/N$ scaling factor appears again in the equality 
\eqref{sigma_fft}, 
due to the way we have chosen the normalization of the discrete Fourier
transform.
This fact is consistent with observation \eqref{w_init}.

The batch normalization defined in \cite{batch-norm}, works with real numbers.
However, we can extend it to complex ones. That is, we may define 
$\layer{B_N}$ be the multivariate function $\layer{B_N} \colon \C{N\times p} \to
\C{N\times p}$ such that $\layer{B_N}(\bfz) = (\layer{B_N}^1(\bfz),\ldots,
\layer{B_N}^p (\bfz))$, $\layer{B_N}^k \colon \C{N\times p} \to \C{N}$ given by:

\begin{equation}
\layer{B_N}^k(\bfz) = \gamma_k \frac{(\bfz^k - \mu_k)}{\sqrt{\sigma_k^2 +
\epsilon}} + \beta_k
\end{equation}
where $\gamma_k \in \C{\space}$ and $\beta_k \in \CN$ are some new weights and
$\epsilon$ is a constant added for numerical stability.

Using again the DFT linearity and the equality \eqref{sigma_fft}, we can show
that:

\begin{equation*}
\begin{aligned}
\FN(\layer{B_N}^k(\bfz)) &= \gamma_k 
\frac{(\FN(\bfz^k) - \FN(\mu_k))}
   {\sqrt{{\frac{1}{N} \sigma_k^2}^\prime + \epsilon}} + \FN(\beta_k)
\\
&= \sqrt{N} \gamma_k \frac{(\FN(\bfz^k) - \mu_k^{\prime})}
   {\sqrt{{\sigma_k^2}^\prime + \epsilon^\prime}} + \FN(\beta_k).
\end{aligned}
\end{equation*}

That is, the batch normalization layer is commuting with the DFT.

\subsection{The Residual Layer Commutes with DFT} \label{residual}

One can show that the residual layer $\layer{R_N}$ as introduced in
\cite{residual}, is commuting with the DFT when modified to use circular
convolutions. That is, if 
\begin{equation}
\layer{R_N}(\bfx, \bfw_1, \bfw_2)
=\conv{\layer{A_N}(\conv{\bfx}{\bfw_1})}{\bfw_2} +
\bfx
\end{equation}
 then
\begin{equation*}
\begin{aligned}
\FN(\layer{R_N}) =& \FN(\conv{\layer{A_N}(\conv{\bfx}{\bfw_1})}{\bfw_2} + \bfx) 
\\ =& \layer{A_{\brac{N}}}(\FN(\bfx) \cdot \FN(\bfw_1))\cdot \FN(\bfw_2) 
\\ &+ \FN(\bfx).
\end{aligned}
\end{equation*}

\subsection{On Complex Weights Initialization} \label{sec:weight-init}

The usual way to initialize the weights in most deep learning frameworks is to
sample them from the (multinomial) normal distribution with zero mean and a
small standard deviation, see \cite{lecun-init}. As we would like to initialize
the weights of Hadamard layers in an equivalent way, we need to understand how 
is this distribution changing under the discrete Fourier transforms. More
precisely, if $\bfx = (x_0,x_1,\ldots, x_{N-1})\in\R{N}$ is a white Gaussian
noise signal of standard deviation $\sigma$, that is, $x_j \sim
\mathcal{N}(0,\,\sigma^{2})$ then the question is what is the distribution of
 $X_k = \FN(x)_k.$ If we note by $Re(z)$
and $Im(z)$ the real and respective the imaginary part of a complex number $z$
then one can show that both $Re(X_i)$ and $ Im(X_i)$ are normally
distributed with standard deviation $\sigma \sqrt{N/2}$, that is, $Re(X_i)
\text{ and } Im(X_i) \sim \mathcal{N}(0,\,\sigma^{2} N/2).$ This results
suggests that we should initialize the complex weights $w_j=u_j +i v_j$
such that 
\begin{equation} \label{w_init}
u_j \sim \mathcal{N}(0,\, 2\sigma^{2}/N) \text{ and } v_j \sim \mathcal{N}(0,\,
2\sigma^{2}/N).
\end{equation}
Note that if the real parts $u_j$ and the imaginary parts $v_j$ are independent 
then $w_j$ are circularly symmetric  complex Gaussians as $u_j$ and $v_j$
have the same variance, see \cite{complex-circ}.

Note also that the term $N/2$ is due to the way we have chosen the normalization
of the discrete Fourier transform. In general, the Fourier transform of the
corresponding probability density of $f(x\mid 0, \sigma^2) = \frac{1}{\sqrt{2\pi
\sigma^2}}e^{\frac{-x^2}{2\sigma^2}}$ is a real Gaussian function that has width
$\sigma$ inverted in the frequency domain:
\begin{equation*}
\begin{aligned}
\mathcal{F}(f(x\mid 0, \sigma^2))(y) &= \int_{-\infty}^{\infty} f(x\mid 0,
\sigma^2) e^{- i y x} dx
\\
&= e^\frac{-\sigma^2y^{2}}{2}.
\end{aligned}
\end{equation*}

\end{document}